\documentstyle[12pt]{amsart}

\newtheorem{theorem}{Theorem}[section]
\newtheorem{lemma}[theorem]{Lemma}

\theoremstyle{definition}
\newtheorem{defn}[theorem]{Definition}

\theoremstyle{remark}
\newtheorem{remark}[theorem]{Remark}

\def\be{\begin{equation}}
\def\ee{\end{equation}}
\def\bea{\begin{eqnarray}}
\def\eea{\end{eqnarray}}

\newcommand{\HE}{\mathrm {HE}}
\newcommand{\Mat}{\mathrm {Mat}}
\newcommand{\HP}{\mathrm {HP}}

\newcommand{\SU}{\mathrm {SU}}
\newcommand{\SO}{\mathrm {SO}}

\newcommand{\GL}{\mathrm {GL}}
\newcommand{\id}{\mathrm {Id}}
\newcommand{\Hom}{\mathrm {Hom}}
\newcommand{\ad}{\mathrm {ad}}
\newcommand{\Tot}{\mathrm {Tot}}

\begin{document}

\begin{center}

\null
\vskip-1truecm
\rightline{IC/97/185}
\vskip1truecm
United Nations Educational Scientific and Cultural Organization\\
and\\
International Atomic Energy Agency\\
\medskip
INTERNATIONAL CENTRE FOR THEORETICAL PHYSICS\\
\vskip1.5truecm
{\bf NON-COMMUTATIVE CHERN CHARACTERS\\ OF COMPACT LIE GROUP
C*-ALGEBRAS} \\
\vskip1.5truecm
Do Ngoc Diep\\
Institute of Mathematics, National Centre for Natural Science and Technology,\\
P.O. Box 631, Bo Ho, 10.000,
Hanoi, Vietnam,
\footnote{\normalsize E-mail: dndiep@@thevinh.ncst.ac.vn}\\
\bigskip
Aderemi O. Kuku\\
International Centre for Theoretical Physics,\\ ICTP P. O. Box 586, 34100,
Trieste, Italy
\footnote{\normalsize E-mail: kuku@@ictp.trieste.it}\\
\bigskip
and\\
\bigskip
Nguyen Quoc Tho\\
Department of Mathematics, Vinh University, Vinh City, Vietnam,\\
and c/o: \\
Institute of Mathematics, National Centre for Natural Science and Technology,\\
P.O. Box 631, Bo Ho, 10.000, Hanoi,
Vietnam.\footnote{\normalsize E-mail: dndiep@@thevinh.ncst.ac.vn}
\end{center}
\vskip0.5truecm
\centerline{ABSTRACT}
\bigskip

For compact Lie groups the Chern characters $ch: K^*(G) \otimes {\mathbf
Q} \to
H^*_{DR}(G;{\mathbf Q})$ were already constructed. In this paper we
propose the
corresponding non-commutative Chern characters, which are also
homomorphisms from
quantum $K$-groups into entire current periodic cyclic homology
groups of group C*-algebras
$ch_{C^*} : K_*(C^*(G)) \to HE_*(C^*(G)).$ We obtain also the
corresponding
algebraic version $ch_{alg} : K_*(C^*(G)) \to HP_*(C^*(G)),$ which
coincides with
the Fedosov-Cuntz-Quillen formula for Chern characters.
\vskip1.5truecm
\begin{center}
MIRAMARE -- TRIESTE\\
\medskip
November 1997
\end{center}

\newpage
\baselineskip=18pt

\section{Introduction}
Let $G$ be a compact Lie group, $H^*_{DR}(G;{\mathbf Q})$ the rational
de Rham
cohomology. It is classical that the Chern character $$ch: K^*(G) \otimes
{\mathbf Q}
\to H^*_{DR}(G;{\mathbf Q}),$$ is an isomorphism, see \cite{W}.There have
also
been other Chern characters from K-theory to such (co)homologies as
Hochschild and
cyclic (co)homologies, etc and the theories have played important role in
geometry, index
theory etc. More recently, some of the theories e.g. cyclic cohomology,
have been
quantized as the so called non-commutative geometry. Since then, there has
been more need
to define non-commutative version of Chern characters that would play in
quantum theories
similar roles to the ordinary commutative Chern characters in classical
theories. The aim of
this paper is to construct such non-commutative Chern characters. In \S2 we
review some
preliminaries on $C^*(G)$ for compact Lie groups $G$ and $KK$ functors. We
also
define the
periodic cyclic homology $\HP_*(A)$ and entire cyclic homology $\HE_*(A)$
for
a general
involutive Banach algebra $A$. 

We define an entire cyclic homology $\HE_*(A)$ of an involutive Banach
algebra $A$ via a family of ideals with ad-invariant trace and via
projective and inductive topologies to form non-commutative analogues of
de Rham currents. We use the dual to the entire cochain of A. Connes at
the level of ideals with ad-invariant trace. When $A=C^*(G)$, where $G$ is
a compact Lie group, our theory coincides with the usual de Rham
cohomology yielding the calculation in the work of Watanabe. (See
\cite{W}, \cite{W'}). Some of these calculations are briefly reviewed in
3.3.

In \S3, we define a non-commutative
Chern
character $ch: K_*(A) \to \HE_*(A)$ (see 3.1) and show in 3.2 that when
$A=C^*(G)$
for a compact Lie group $G$, $ch$ is an isomorphism, which can be
identified
with the
classical Chern character $K_*^W({\mathbf C}({\mathbf T})) \to
\HE^W_*({\mathbf
C}({\mathbf T}))$ that is also an isomorphism where ${\mathbf T}$ is a
maximal torus
of $G$ with
Weyl group $W$. In \S4, we indicate at first that the periodic cyclic
homology
$\HP_*(A)$
defined in \S2 through ideals of $A$ coincides with the periodic cyclic
homology of
Cuntz-Quillen \cite{CQ} when $A=C^*(G)$ for a compact Lie group $G$. We
then go
on to show in
4.5 that the Chern character $ch: K_*(C^*(G)) \to \HP_*(C^*(G))$ is an
isomorphism
which can be identified with the classical Chern character $K_*({\mathbf
C}({\mathbf T}))
\to \HP_*({\mathbf C}({\mathbf T}))$ which is also an isomorphism.

{\it Notes on Notation}:
For any compact Lie group $G$, we write $K^*(G)$ for the ${\mathbf Z}/2$
graded topological K-theory of $G$. We use Swan's theorem to identify
$K^*(G)$ with the ${\mathbf Z}/2$-graded $K_*(C(G))$. For any involutive
Banach algebra $A$, $K_*(A)$, $\HE_*(A)$, $\HP_*(A)$, are 
${\mathbf Z}/2$-graded algebraic K-group of $A$, entire cyclic homology of
$A$, and periodic cyclic homology of $A$, respectively. If ${\mathbf T}$
is a maximal torus of $G$ with the Weyl group $W$, write ${\mathbf
C}({\mathbf T})$ for the ring of complex valued functions on ${\mathbf
T}$. We write $K_*^W({\mathbf C}({\mathbf T}))$, (resp. $\HE_*^W({\mathbf
C}({\mathbf T}))$ for the ${\mathbf Z}/2$-graded $W$-equivariant K-group,
(resp. entire cyclic homology) of ${\mathbf C}({\mathbf T})$. 

\section{Preliminaries}

\subsection{Group C*-Algebras} We introduce in this subsection the
definition and the
main properties of C*-algebras of compact Lie groups we need for later
use.
Let us recall that in a compact Lie group $G$ there exists a (bi)invariant
Haar
measure denoted
$dg$. Consider the involutive Banach algebra $L^1(G)$ of absolutely
integrable
complex-valued
functions on $G$ with the convolution product $$(f*g)(x) := \int_G
f(y)g(y^{-1}x)dy,$$
involution
$$f^*(x) := \overline{f(x^{-1})}$$
and the $L^1$-norm
$$\Vert f\Vert := \int_G \vert f(x)\vert dx.$$ This algebra is not regular
in the sense that in
general $$\Vert \varphi * \psi \Vert \ne \Vert \varphi \Vert \Vert \psi
\Vert.$$ Let us consider
the new regular norm defined by $$\Vert \varphi \Vert := \sup_{\pi\in
\hat{G}}\Vert
\pi(\varphi)\Vert,$$ where by definition, $\hat{G}$ is the dual object of
$G$, i.e. the set of
(unitary) equivalence classes of irreducible (unitary) representations of
$G$. The
completion of the
involutive Banach algebra $L^1(G)$ with respect to this norm is a regular
complete involutive
Banach algebra with unit element, which is called the {\it group
C*-algebra} of
$G$.

The following facts are well known, see \cite{K}.
\begin{enumerate}
\item The dual object $\hat{G}$ is not more than denumerable, finite if the
group $G$ is
finite and of cardinal ${\mathbf N}$ if $G$ is compact.
\item Every irreducible representation of the compact Lie group $G$ is
equivalent to a unitary
one (unitarization) and is finite dimensional, of dimension $n_i$, say
\item Every representation $\pi$ of $G$ can be extended to a
*-representation of
$C^*(G)$ by
the Fourier-Gel'fand transform $$\hat{f}(\pi) = \pi(f) := \int_G \pi(x)f(x)
dx.$$ \item There
is one-to-one correspondence between the irreducible representations
$\pi$ of $G$
and the
non-degenerate irreducible *-representation of $C^*(G)$.
\item Let us fix a bijection between $\pi_n\in\hat{G}$ and $n\in{\mathbf
N}$. Then, there exists a constant $c=c_f$, where $f$ is in $C^*(G)$, such
that $$\Vert
\hat{f}(n) -
c_f.\id\Vert \to 0,$$ as $n \to \infty$. \end{enumerate}
Let us denote the {\it restricted} Cartesian product of matrix algebras
$\Mat_{n_i}({\mathbf C})$ by
$${\prod'}_{i=1}^\infty\Mat_{n_i}({\mathbf C})= \{\hat{f}\quad;\quad \Vert
\hat{f}(n) -
c_f.\id\Vert \to 0, \mbox{ as }n \to \infty \}$$ The main property of the
group C*-algebra
$C^*(G)$ is the fact that $$C^*(G) \cong
{\prod'}_{i=1}^\infty\Mat_{n_i}({\mathbf
C}).$$ This means also that
$$C^*(G) = \varinjlim I_N,$$ and that the right-hand side is dense in norm,
where $$I_N
:= {\prod}_{i=1}^N\Mat_{n_i}({\mathbf C})$$ are two-sided closed ideals in
$C^*(G)$.

\subsection{KK-groups} We recall in this subsection a brief definition of
operator
K-functors.  Following G. G. Kasparov, a
Fredholm representation of a C*-algebra $A$ is a triple
$(\pi_1,\pi_2,F)$,
consisting of *-representations $\pi_1, \pi_2 : A \to {\mathcal
L}({\mathcal H}_B)$ and a
Fredholm operator $F\in {\mathcal F}({\mathcal H}_B)$,
admitting an adjoint operator, on the
Hilbert
C*-module ${\mathcal H}_B= \ell^2_B$ over C*-algebra $B$, satisfying the
relations
$$\pi_1(a)F -
F\pi_2(a) \in {\mathcal K}_B,$$ where ${\mathcal K}_B$ is the ideal of compact
(adjointable) C*-module endomorphisms of ${\mathcal H}_B$.
The classes of homotopy invariance and unitary equivalence of Fredholm
modules form the
so called Kasparov
operator K-group $KK^*(A,B)$. Herewith put $A=C^*(G)$, $B =
{\mathbf C}$,
where $G$ is a compact Lie group. In this case we have $K_*(C^*(G)) =
KK^*(C^*(G),{\mathbf C})$, where $K_*(A)$ is algebraic $K$-group of $A$.
Note
that
$$KK_0(A, {\mathbf C}) \cong K_0(A),$$
$$KK_1(A, {\mathbf C}) \cong K_1(A).$$
\subsection{Entire Homology of Involutive Banach Algebras}Let $A$ be an
involutive
Banach algebra. Recall that A. Connes defines entire cyclic
cohomology $\HE^*(A)$ and a pairing
$$K_*(A) \times \HE^*(A) \to {\mathbf C}.$$ Also M. Khalkhali \cite{Kh1}
proved
Morita and homotopy invariance of $\HE^*(A)$. We now define the entire
homology
$\HE_*(A)$ as follows: Given a collection $\{I_\alpha\}_{\alpha\in\Gamma}$
of ideals in $A$, equipped with a so called $\ad_A$-invariant trace
$$\tau_\alpha : I_\alpha
\to {\mathbf C},$$ satisfying the properties: \begin{enumerate}
\item $\tau_\alpha$ is a {\it continuous linear} functional, normalized as
$\Vert\tau_\alpha\Vert = 1$,
\item $\tau_\alpha$ is {\it positive} in the sense that $$\tau_\alpha(a^*a)
\geq 0, \forall
\alpha \in \Gamma,$$ where the map $a \mapsto a^*$ is the involution
defining the
involutive Banach algebra structure, i.e. an anti-hermitian endomorphism
such that $a^{**}
= a$ \item $\tau_\alpha$ is {\it strictly positive} in the sense that
$\tau_\alpha(a^*a) = 0$, iff
$a=0$, for every $\alpha\in \Gamma$. \item $\tau_\alpha$ is $\ad_A$-{\it
invariant} in the
sense that $$\tau_\alpha(xa) = \tau_\alpha(ax), \forall x\in A, 
a\in I_\alpha,$$
\end{enumerate}
then we have for every $\alpha \in \Gamma$ a scalar product $$\langle
a,b\rangle_\alpha
:=\tau_\alpha(a^*b)$$ and also an inverse system $\{I_\alpha,
\tau_\alpha\}_{\alpha \in
\Gamma}$. Let $\bar{I}_\alpha$ be the completion of $I_\alpha$ under the
scalar product
above and $\widetilde{\bar{I}_\alpha}$ denote $I_\alpha$ with formally
adjoined unity
element. Define the $C^n(\widetilde{\bar{I}_\alpha})$ the set of
$n+1$-linear maps $\varphi :
(\widetilde{\bar{I}_\alpha})^{n+1} \to {\mathbf C}$.

For $\alpha \leq \beta$, we have a well-defined map $$D^\beta_\alpha :
C^n(\widetilde{\bar{I}_\alpha}) \to C^n(\widetilde{\bar{I}_\beta}),$$ which
makes $\{
C^n(\widetilde{\bar{I}}_\alpha)\}$ into a direct system. Write $Q = \varinjlim
C^n(\widetilde{\bar{I}_\alpha})$. Remark that it admits a Hilbert space 
structure, see \cite{DT1}-\cite{DT2}. Let $C_n(A) := \Hom(\varinjlim
C^n(\widetilde{\bar{I}_\alpha}), {\mathbf C})$.

Let $$b, b' : C^n(\widetilde{\bar{I}}_\alpha) \to
C^{n+1}(\widetilde{\bar{I}}_\alpha),$$
$$N : C^n(\widetilde{\bar{I}}_\alpha) \to
C^n(\widetilde{\bar{I}}_\alpha),$$ $$\lambda :
C^n(\widetilde{\bar{I}}_\alpha) \to C^n(\widetilde{\bar{I}}_\alpha),$$
$$S : C^{n+1}(\widetilde{\bar{I}}_\alpha)\to
C^n(\widetilde{\bar{I}}_\alpha)$$ be defined as in A. Connes \cite{Ca}.
We follow the notations in \cite{Kh1} . Denotes by $b^*, (b')^*, N^*,
\lambda^*, S^*$ the corresponding adjoint operators.

We now have a bi-complex
$$\begin{array}{ccccccccc}
& &\phantom{(-b')^*}\vdots & &\phantom{b^*}\vdots & &\phantom{(-b')^*}\vdots & & \\
& &(-b')^*\downarrow & &b^*\downarrow & &(-b')^*\downarrow & & \\
& 1-\lambda^*& & N^*& & 1-\lambda^*& &N^* & \\ \dots& \longleftarrow & \phantom{(-
b')^*}C_1(A) & \longleftarrow & \phantom{b^*}C_1(A) &\longleftarrow & \phantom{(-
b')^*}C_1(A) & \longleftarrow & \ldots\\
& & & & & & & & \\
& &(-b')^*\downarrow & &b^*\downarrow & &(-b')^*\downarrow & & \\
& 1-\lambda^*& & N^*& & 1-\lambda^*& &N^* & \\ \dots& \longleftarrow & \phantom{(-
b')^*}C_0(A) & \longleftarrow & \phantom{b^*}C_0(A) &\longleftarrow & \phantom{(-
b')^*}C_0(A) & \longleftarrow & \ldots\\
\end{array}\leqno{{\mathcal C}(A):}$$
with $b^*$ in the even columns and $(-b')^*$ in the odd columns, where *
means the corresponding adjoint operator. Now we have
$$\Tot({\mathcal C}(A)^{even}) = \Tot({\mathcal C}(A)^{odd}) :=
\oplus_{n\geq 0}
C_n(A),$$ which is periodic with period two. Hence, we have
$$
\bigoplus_{n\geq 0} C_n(A)\begin{array}{c}\partial\\ \longleftarrow\\
\longrightarrow\\ \partial \end{array}\bigoplus_{n\geq 0} C_n(A), $$
where $\partial = d_v + d_h$ is the total differential.

Let $\HP_*(A)$ be the homology of the total complex $(\Tot{\mathcal C}(A))$.
Remark that this $\HP_*(A)$ is, in general different from the 
$\HP_*(A)$ of Cuntz-Quillen. 

\begin{defn}
An even (or odd) chain $(f_n)_{n\geq 0}$ in ${\mathcal C}(A)$ is called
{\it entire} if the
radius of convergence of the power series $\sum_n \frac{n!}{[\frac{n}{2}]!}
\Vert
f_n\Vert z^n$, $z\in {\mathbf C}$ is infinite.
\end{defn}

Let $C_e(A)$ be the subcomplex of $C(A)$ consisting of entire chains. Then
we have a
periodic complex.

\begin{theorem} Let
$$\Tot(C_e(A)^{even}) = \Tot(C_e(A)^{odd}) := \bigoplus_{n\geq 0} C^e_n(A),$$
where $C^e_n(A)$ is the entire $n$-chain. Then we have a complex of
entire chains with the total differential $\partial $
$$
\bigoplus_{n\geq 0} C^e_n(A) \begin{array}{c}\partial\\ \longleftarrow\\
\longrightarrow\\ \partial \end{array} \bigoplus_{n\geq 0}
C^e_n(A)$$
\end{theorem}
The homology of this complex is called also the {\it entire homology} and
denoted by $\HE_*(A)$. Note that this entire homology is defined through
the inductive limits of ideals with ad-invariance trace. 

\section{Non-commutative Chern Characters for Involutive Banach Algebras}

Let $A$ be an involutive Banach algebra. In this section, we construct a
non-commutative
character
$$ch_{C^*} : K_*(A) \to \HE_*(A)$$ and later show that when $A= C^*(G)$,
this Chern character reduces up to isomorphism to classical Chern
character.

Let $A$ be an involutive Banach algebra with unity. 
\begin{theorem}
There exists a Chern character $$ch_{C^*} : K_*(A) \to \HE_*(A).$$ 
\end{theorem}
\begin{pf}
We first recall that there exists a pairing $$K_n(A) \times C^n(A) \to
{\mathbf C}$$ due to
A. Connes, see \cite{Co}. Hence there exists a map
$ K_n(A) \stackrel {C_n }{ \longrightarrow }\Hom(C^n(A),{\mathbf C}).$ So,
by 1.1,
we have for each $\alpha \in {\mathbf C}$, a map $K_n(A)
\stackrel{C_n^\alpha}{\longrightarrow}
\Hom(C^n(\widetilde{\bar{I}_\alpha}),{\mathbf C})$ and hence a map $K_n(A)
\stackrel{C_n}{\longrightarrow}
\Hom(\varinjlim_\alpha C^n(\widetilde{\bar{I}_\alpha}),{\mathbf C}).$ We
now show that $C_n$
induces
the Chern map
$$ch: K_n(A) \to \HE_n(A)$$

Now let $e$ be an idempotent in $M_k(A)$ for some $k\in {\mathbf N}$. It
suffices to
show that for $n$ even, if $\varphi = \partial\psi$, where $\varphi\in
C^n(\widetilde{\bar{I}_\alpha})$ and $\psi \in
C^{n+1}(\widetilde{\bar{I}_\alpha})$,
then $$\langle e,\varphi\rangle = \sum_{n=1}^\infty
\frac{(-1)^n}{n!}\varphi(e,e,\dots,e)
= 0.$$
However, this follows from Connes' results in (\cite{Co},Lemma 7).

The proof of the case for $n$ odd would also follow from \cite{Co}. \end{pf}

Our next result computes the Chern character in 2.1 for $A=C^*(G)$ by
reducing it to the
classical case.

\begin{theorem}
Let ${\mathbf T}$ be a fixed maximal torus of $G$ with Weyl group $W:=
N_G({\mathbf T})/{\mathbf T}$.
Then the Chern character $$ch_{C^*} : K_*({\mathbf C}^*(G)) \to
\HE_*(C^*(G))$$ is an
isomorphism, which can be identified with the classical Chern character
$$ch: K^W_*({\mathbf C}({\mathbf T})) \to \HE^W_*({\mathbf C}({\mathbf
T}))$$ that
is also an isomorphism.
\end{theorem}
\begin{pf}
First observe that we have an isomorphism $$K_*(C^*(G)) \cong 
K_*({\prod'}_{i=1}^\infty\Mat_{n_i}({\mathbf C})
= \varinjlim K_*(\prod_{i=1}^N \Mat_{n_i}({\mathbf C})) \cong $$
$$\cong K_*(\prod_{\lambda=\mbox{highest weight mod W}} {\mathbf
C}_\lambda) =
K^W_*({\mathbf
C}({\mathbf T})).$$
Next we have
$$\HE_*(C^*(G)) \cong \HE_*(\varinjlim\prod_{i=1}^N \Mat_{n_i}({\mathbf
C})) =
\varinjlim\HE_*(\prod_{i=1}^N \Mat_{n_i}({\mathbf C}))$$ $$\cong
\HE_*(\prod_{\lambda =\mbox{highest weight mod W}} {\mathbf C}_\lambda)
=\HE^W_*({\mathbf
C}({\mathbf T})),$$
where ${\mathbf 
C_\lambda} = {\mathbf C}$ are enumerated by the highest weights of the
corresponding irreducible representations of $G$.

Because the irreducible representations of compact Lie groups are defined
by their characters (highest weights), the above product is indexed by
 highest weights modulo the action of the Weyl group $W$.

Furthermore, we have from standard results in topology that
$$K_W^*({\mathbf T}) \cong K^*_W(B{\mathbf T}).$$ 

Now consider the commutative diagram
$$\begin{array}{ccc}
K_*(C^*(G)) & \stackrel{\eta}{----\longrightarrow} & K_*^W({\mathbf
C}({\mathbf
T}))\\
\vert\phantom{ch_{C^*}} & &\vert\phantom{ch}\\ \vert ch_{C^*} &
&\vert ch\\ \downarrow \phantom{ch_{C^*}} & & \downarrow \phantom{ch} \\
\HE_*(C^*(G)) &
\stackrel{\delta}{----\longrightarrow} & \HE_*^W({\mathbf C}({\mathbf T}))
\end{array}$$
Since $\eta$, $\delta$ are isomorphisms and $ch$ is also an isomorphism, we
have that
$ch_{C^*} = \delta^{-1}(ch)\eta$ is also an isomorphism. \end{pf}

\begin{remark}
For some classical groups, e. g. $\SU(n+1)$, $\SO(2n+1)$, $\SU(2n)$,
$Sp(n)$ etc. the
groups $$K^*(G) \cong K_*^W({\mathbf C}({\mathbf T}))\cong  K^*_W({\mathbf
T}) \cong K_*(C^*(G))\cong$$
$$\cong \HE_*(C^*(G))\cong \HE^W_*({\mathbf C}({\mathbf T})\cong
H^W_*({\mathbf C}({\mathbf T})\cong H^*_W({\mathbf
T}) \cong H^*(G) \cong \HP_*(C^*(G))$$ are
as follows.

(a). For any compact Lie group $G$, let $R[G]$ be the representation ring.
Then
the
${\mathbf Z}/(2)$-graded algebra $K^*(G) =
\wedge_{\mathbf
C}(\beta(\rho_1), \beta(\rho_2),\dots,\beta(\rho_n))$, where $\rho_i$ are
the standard
irreducible representations and $\beta :R[G] \to K^*(G)$ is the Bott map.
Hence, from
\cite{W} we have $$K^*(\SU(n+1)) \cong \wedge_{\mathbf C}(\beta(\rho_1),
\dots,
\beta(\rho_n)),$$ $$K^*(\SO(2n+1)) \cong \wedge_{\mathbf C}(\beta(\rho_1),
\dots,
\beta(\rho_n), \varepsilon_{2n+1}).$$

(b). It follows from \cite{W'} that the ${\mathbf Z}/(2)$-graded complex
cohomology
groups are exterior algebras over ${\mathbf C}$ and in particular
$$H^*(\SU(2n)) \cong \wedge_{\mathbf C}(x_3,x_5,\dots, x_{4n-1}),$$
$$H^*(Sp(n))
\cong \wedge_{\mathbf C}(x_3,x_7,\dots, x_{4n-1}),$$ $$H^*(\SU(2n+1)) \cong
\wedge_{\mathbf C}(x_3,x_5,\dots, x_{4n+1}),$$ $$H^*(\SO(2n+1)) \cong
\wedge_{\mathbf C}(x_3,x_7,\dots, x_{4n-1}).$$

(c). Define a function $\Phi: {\mathbf N} \times {\mathbf N} \times
{\mathbf N} \to
{\mathbf Z}$ by
$$\Phi(n,k,\varepsilon) = \sum_{i=1}^k (-1)^{i-1} \binom{n}{k-i} i^{q-1}.$$
It then
follows from \cite{W}, \cite{W'} that we have Chern character $ch:
K^*(\SU(n+1)) \to H^*(SU(n+1))$ given by
$$ch(\beta(\rho_k)) =
\sum_{i=1}^n \frac{(-1)^i}{i!}\Phi(n+1,k,i+1) x_{2i+1}, \forall k \geq 1,$$
$ch : K^*(\SO(2n+1)) \to
H^*(\SO(2n+1)),$ given by the formula $$ch(\beta(\lambda_k)) =
\sum_{i=1}^n
\frac{(-1)^{i-1}2}{(2i-1)!}\Phi(2n+1,k,2i)x_{4i-1} (\forall k=1,2,\dots,n-1)$$
$$ch(\varepsilon_{2n+1}) = \sum_{i=1}^n
\frac{(-1)^{i-1}}{2^{n-1}(2i-1)!}\sum_{k=1}^n\Phi(2n+1,k,2i)x_{4i-1}.$$
\end{remark}

\section{Algebraic Version of Non-commutative Chern Character}

Let $G$ be a compact Lie group, $\HP_*(C^*(G))$ the periodic cyclic
homology
introduced
in \S2. Since $C^*(G) = \varinjlim_N \prod_{i=1}^N \Mat_{n_i}({\mathbf C})$,
$\HP_*(C^*(G))$ coincides with the $\HP_*(C^*(G))$ defined by J. Cuntz-D.
Quillen
\cite{CQ}.

\begin{lemma}\label{31}
Let $\{I_N\}_{N\in {\mathbf N}}$ be the above defined collection of ideals
in $C^*(G)$.
Then $$K_*(C^*(G)) = \varinjlim_{N \in {\mathbf N}} K_*(I_N) =
K^W_*({\mathbf
C}({\mathbf T})),$$ where ${\mathbf T}$ is the fixed maximal torus of $G$.
\end{lemma}
\begin{pf}
First note that the algebraic K-theory of $C^*$-algebras has the stability
property
$$K_*(A \otimes M_n({\mathbf C})) \cong K_*({\mathbf C}).$$
Hence,
$$\varinjlim K_*(I_{n_i}) \cong K_*(\prod_{\lambda =\mbox{highest weight
mod W}}
{\mathbf C}_\lambda) \cong
K^W_*({\mathbf C}({\mathbf T}))$$ by Pontryagin duality, where ${\mathbf 
C_\lambda} = {\mathbf C}$ are enumerated by the highest weights of the
corresponding irreducible representations of $G$.

\end{pf}

J. Cuntz and D. Quillen \cite{CQ} defined the so called $X$-complexes of
${\mathbf C}$-algebras and then used some ideas of Fedosov product to
define algebraic Chern
characters. We now briefly recall their definitions. For a
(non-commutative) associate ${\mathbf C}$-algebra $A$, consider the space
of even
non-commutative differential forms $\Omega^+(A) \cong RA$, equipped with
the Fedosov
product
$$\omega_1 \circ \omega_2 := \omega_1\omega_2 - (-1)^{\vert \omega_1\vert}
d\omega_1
d\omega_2,$$ see \cite{CQ}. Consider also the ideal $IA := \oplus_{k\geq 1}
\Omega^{2k}(A)$. It is easy to see that $RA/IA \cong A$ and that $RA$
admits the
universal property that any based linear map $\rho : A \to M$ can be
uniquely extended to a
derivation $D : RA \to M$. The derivations $D : RA \to M$ bijectively
correspond to lifting
homomorphisms from $RA$ to the semi-direct product $RA \oplus M$, which also
bijectively correspond to linear map $\bar\rho : \bar{A}= A /{\mathbf C}
\to M$ given by $$
a\in \bar{A} \mapsto D(\rho a).$$ From the universal property of
$\Omega^1(RA)$, we
obtain a bimodule isomorphism $$RA \otimes \bar{A} \otimes RA \cong
\Omega^1(RA).$$ As in \cite{CQ}, let $\Omega^-A = \oplus_{k\geq 0}
\Omega^{2k+1}A$. Then we have $$\Omega^{-}A \cong RA \otimes \bar{A} \cong
\Omega^1(RA)_\# := \Omega^1(RA)/[(\Omega^1(RA),RA)].$$
\par
J. Cuntz and D. Quillen proved
\begin{theorem}(\cite{CQ}, Theorem1):
There exists an isomorphism of
${\mathbf Z}/(2)$-graded complexes
$$\Phi : \Omega A = \Omega^+A \oplus \Omega^{-}A \cong RA \oplus
\Omega^1(RA)_\#,$$ such that
$$\Phi : \Omega^+A \cong RA,$$ is defined by $$\Phi(a_0da_1\dots da_{2n} =
\rho(a_1)\omega(a_1,a_2) \dots \omega(a_{2n-1},a_{2n}),$$
and $$ \Phi : \Omega^{-}A \cong \Omega^1(RA)_\#,$$ $$\Phi(a_0da_1\dots
da_{2n+1})
= \rho(a_1)\omega(a_1,a_2)\dots \omega(a_{2n-1},a_{2n})\delta(a_{2n+1}).$$
With
respect to this identification, the product in $RA$ is just the Fedosov
product on even
differential forms and the differentials on the $X$-complex
$$X(RA) : \qquad RA\cong \Omega^+A \to \Omega^1(RA)_\# \cong \Omega^{-}A
\to RA
$$ become the operators
$$\beta = b - (1+\kappa)d : \Omega^{-}A \to \Omega^+A,$$ $$\delta =
-N_{\kappa^2} b
+ B : \Omega^+A \to \Omega^{-}A,$$ where $N_{\kappa^2} =
\sum_{j=0}^{n-1} \kappa^{2j}$, $\kappa(da_1\dots da_n) := da_n\dots da_1$.
\end{theorem}
\par
Let us denote by $IA \triangleleft RA$ the ideal of even non-commutative
differential forms
of
order $\geq 2$. By the universal property of $\Omega^1$ $$\Omega^1(RA/IA) =
\Omega^1RA/((IA)\Omega^1RA + \Omega^1RA.(IA) + dIA).$$ Since $\Omega^1RA =
(RA)dRA = dRA.(RA)$, then $\Omega^1RA(IA) \cong IA\Omega^1RA \mod
[RA,\Omega^1R].$
$$\Omega^1(RA/IA)_\# = \Omega^1RA /([RA,\Omega^1RA]+IA.dRA + dIA).$$ For
$IA$-adic tower $RA/(IA)^{n+1}$, we have the complex $${\mathcal
X}(RA/(IA)^{n+1}) : \qquad RA/IA^{n+1} \leftarrow
\Omega^1RA/([RA,\Omega^1RA]+(IA)^{n+1}dRA + d(IA)^{n+1}).$$
Define
$${\mathcal X}^{2n+1}(RA,IA) : \quad RA/(IA)^{n+1} \to
\Omega^1RA/([RA,\Omega^1RA]+(IA)^{n+1}dRA + d(IA)^{n+1}) $$ $$\to
RA/(IA)^{n+1},$$
$${\mathcal X}^{2n}(RA,IA): \quad RA/((IA)^{n+1} +[RA,IA^n]) \to
\Omega^1RA/([RA,\Omega^1RA]+d(IA)^ndRA)$$ $$\to RA/((IA)^{n+1}
+[RA,IA^n]).$$
One has
$$b((IA)^ndIA) = [(IA)^n,IA] \subset (IA)^{n+1},$$ $$d(IA)^{n+1} \subset
\sum_{j=0}^n (IA)^jd(IA)(IA)^{n-j} \subset (IA)^n dIA
+ [RA,\Omega^1RA].$$
and hence
$${\mathcal X}^1(RA,IA) = X(RA,IA),$$
$${\mathcal X}^0(RA,IA) = (RA/IA)_\#.$$
There is a sequence of maps between complexes $$\dots \to X(RA/IA) \to
{\mathcal
X}^{2n+1}(RA,IA)\to {\mathcal X}^{2n}(RA,IA) \to X(RA/IA) \to \dots $$ We
have the
inverse limits
$$\hat{X}(RA,IA) := \varprojlim X(RA/(IA)^{n+1}) = \varprojlim {\mathcal
X}^n(RA,IA).$$
Remark that
$${\mathcal X}^q = \Omega A/F^q\Omega A,$$ $$\hat{X}(RA/IA)= \hat{\Omega}A.$$

We quote the second main result of J. Cuntz and D. Quillen (\cite{CQ},
Thm2), namely:

$$H_i\hat{\mathcal X}(RA,IA) = \HP_i(A).$$

We now apply  this machinery to our case. First we have the following.
\begin{lemma}
$$\varinjlim \HP^*(I_N) \cong \HP^*({\mathbf C}({\mathbf T})).$$ \end{lemma}
\begin{pf} By similar arguments as in the previous lemma \ref{31}. More
precisely, we
have $$\HP(I_{n_i}) = \HP(\prod_{\lambda =\mbox{highest weight mod W}}
{\mathbf
C}_\lambda)
\cong
\HP^W_*({\mathbf C}({\mathbf T}))$$ by Pontryagin duality.
\end{pf}

Now, for each idempotent $e\in M_n(A)$ there is an unique element $x\in
M_n(\widehat{RA})$.
Then the element $$\tilde{e} := x + (x-\frac{1}{2})\sum_{n\geq 1} \frac{2^n(2n-
1)!!}{n!}(x-x^{2n})^{2n}\in M_n(\widehat{RA})$$ is a lifting of $e$ to an
idempotent
matrix in $M_n(\widehat{RA})$. Then the map $[e] \mapsto tr(\tilde{e})$
defines the map
$K_0(A) \to H_0(X(\widehat{RA})) = \HP_0(A)$. To an element $g\in\GL_n(A)$
one
associates an element $p\in \GL(\widehat{RA})$ and to the element $g^{-1}$
an element
$q\in \GL_n(\widehat{RA})$ then put
$$x = 1- qp, \mbox{ and } y = 1-pq.$$
And finally, to each class
$[g]\in \GL_n(A)$ one associates $$tr(g^{-1}dg) = tr(1-x)^{-1}d(1-x) =
d(tr(log(1-x))) =
-tr\sum_{n=0}^\infty x^ndx\in \Omega^1(A)_\#.$$ Then $[g] \to tr(g^{-1}dg)$
defines
the map $K_1(A) \to HH_1(A) = H_1(X(\widehat{RA})) = \HP_1(A)$.

\begin{defn} Let $\HP(I_{n_i})$ be the periodic cyclic cohomology defined by
Cuntz-Quillen. Then the pairing
$$K_*^{alg}(C^*(G)) \times \bigcup_N \HP^*(I_N) \to {\mathbf C}$$ defines an
algebraic non-commutative Chern character $$ch_{alg} : K_*^{alg}(C^*(G)) \to
\HP_*(C^*(G)),$$ which gives us a variant of non-commutative Chern
characters with
values in $\HP$-groups. \end{defn}

We close this section with an algebraic analogue of theorem 3.2.

\begin{theorem}
Let $G$ be a compact Lie group and ${\mathbf T}$ a fixed maximal
compact torus
of $G$.
Then in the notations of 4.3, the Chern character $$ch_{alg} : K_*(C^*(G)) \to
\HP_*(C^*(G))$$ is an isomorphism, which can be identified with the
classical Chern
character $$ch: K^W_*({\mathbf C}({\mathbf T})) \to \HP^W_*({\mathbf
C}({\mathbf
T}))$$
which is also an isomorphism.
\end{theorem}
\begin{pf}
First note that $$K_*(C^*(G)) \cong K_*(\varinjlim \prod_{i=1}^N
\Mat_{n_i}({\mathbf
C})) = \varinjlim K_*(\prod_{i=1}^N \Mat_{n_i}({\mathbf C})) \cong $$
$$\cong K_*(\prod_{\lambda =\mbox{highest weight mod W}} {\mathbf
C}_\lambda) =
K_*({\mathbf
C}({\mathbf T})).$$
Next we have
$$\HP_*(C^*(G)) \cong \HP_*(\varinjlim\prod_{i=1}^N \Mat_{n_i}({\mathbf
C})) =
\varinjlim\HP_*(\prod_{i=1}^N \Mat_{n_i}({\mathbf C}))$$ $$\cong
\HP(\prod_{\lambda=\mbox{highest weight mod W}} {\mathbf C}_\lambda)
=\HP^W_*({\mathbf
C}({\mathbf
T})).$$
Furthermore, by a result of Cuntz-Quillen for the commutative ${\mathbf
C}$-algebra
$A$, we have a canonical isomorphism from periodic cyclic homology to the
${\mathbf
Z}/(2)$-graded de Rham homology which is an isomorphism when $A$ is smooth.
Hence
$$K^W_*({\mathbf C}({\mathbf T})) \cong \HP^W_*({\mathbf C}({\mathbf
T})).$$
Now we
have a commutative diagram
$$\begin{array}{ccc}
K_*(C^*(G)) & \stackrel{\eta}{----\longrightarrow} & K^W_*({\mathbf
C}({\mathbf T}))\\
\vert\phantom{ch_{alg}} & &\vert\phantom{ch}\\ \vert ch_{alg} &
&\vert ch\\ \downarrow \phantom{ch_{alg}} & & \downarrow
\phantom{ch} \\ \HP_*(C^*(G)) &
\stackrel{\delta}{----\longrightarrow} & \HP^W_*({\mathbf C}({\mathbf T}))
\end{array}$$
where $\eta, \delta, ch$ are isomorphisms. Hence, $ch_{alg} =
\delta^{-1}(ch)\eta$ is an
isomorphism.
\end{pf}

\section*{Acknowledgments}
This work was completed during the stay of the first author as a visiting
mathematician at
the International Centre for Theoretical Physics, Trieste, Italy.
He would like
to thank ICTP for the hospitality, without which this work would not have
been possible.

This work is supported in part by the International Centre for Theoretical
Physics,
Trieste, Italy and the National Foundation for Research in Natural Sciences
of Vietnam.

\end{document}